\documentclass[12pt]{article}
\usepackage{amsmath}
\usepackage{amssymb}
\usepackage{delarray}
\usepackage{amsfonts}
\usepackage{amsthm}
\usepackage[english]{babel}
\usepackage{enumerate}
\usepackage{fancyhdr}

\def \egal {\stackrel{{\rm def}}{=}}

\def\tf{\tilde{f}_n}

\newcommand \cA{{\cal A}}
\newcommand \cB{{\cal B}}
\newcommand \cC{{\cal C}}
\newcommand \cD{{\cal D}}
\newcommand \cF{{\cal F}}

\newcommand \cP{{\cal P}}

\newcommand \cT{{\cal T}}

\newcommand \cX{{\cal X}}

\newcommand \cZ{{\cal Z}}
\newcommand{\1}{{\rm 1}\kern-0.24em{\rm I}}

\theoremstyle{plain}
\newtheorem{thm}{Theorem}
\newtheorem{cor}{Corollary}
\newtheorem{lem}{Lemma}

\newtheorem{rem}{Remark}
\newtheorem{defn}{Definition}

\date{}

\begin{document}

\title{{\bf  Adapting to Unknown Smoothness by Aggregation of Thresholded Wavelet Estimators.}}
\author{
  Christophe Chesneau and Guillaume Lecu\'e\\
{\it Universit\'e Paris VI } \\
}

\maketitle
\date

\begin{abstract}
We study the performances of an adaptive procedure based on a convex
combination, with data-driven weights, of term-by-term thresholded
wavelet estimators. For the bounded regression model, with random
uniform design, and the nonparametric density model, we show that
the resulting estimator is optimal in the minimax sense over all
Besov balls under the $L^2$ risk, without any logarithm factor.
\end{abstract}

\section{Introduction}
Wavelet shrinkage methods have been very successful in nonparametric
function estimation. They provide estimators that are spatially
adaptive and (near) optimal over a wide range of function classes.
Standard approaches are based on the term-by-term thresholds. A
well-known example is the hard thresholded estimator introduced by
\cite{donohoj1}. If we observe $n$ statistical data and if the
unknown function $f$ has an expansion of the form $f=\sum_{j}\sum_k
\beta_{j,k}\psi_{j,k}$ where $\{\psi_{j,k}, \ j,k\}$ is a wavelet
basis and $(\beta_{j,k})_{j,k}$ is the associated wavelet
coefficients, then the term-by-term wavelet thresholded method
consists in three steps:
\begin{enumerate}
\item a linear step corresponding to the estimation of the
coefficients $\beta_{j,k}$ by some estimators $\hat \beta_{j,k}$
constructed from the data, \item a non-linear step consisting in a
thresholded procedure $T_{\lambda}(\hat \beta_{j,k})\1_{\left
\lbrace |\hat \beta_{j,k}|\ge \lambda_j\right \rbrace}$ where
$\lambda=(\lambda_j)_{j}$ is a positive sequence and
$T_{\lambda}(\hat \beta_{j,k})$ denotes a certain transformation of
the $\hat \beta_{j,k}$ which may depend on $\lambda$, \item a
reconstruction step of the form $\hat f_{\lambda}=\sum_{j\in
\Omega_n}\sum_{k}T_{\lambda}(\hat \beta_{j,k})\1_{\left \lbrace
|\hat \beta_{j,k}|\ge \lambda_j\right \rbrace}\psi_{j,k} $ where
$\Omega_n$ is a finite set of integers depending on the number $n$
of data.
\end{enumerate}
Naturally, the performances of $\hat f_{\lambda}$ strongly depend on
the choice of the threshold $\lambda$. For the standard statistical
models (regression, density,...), the most common choice is the
universal threshold introduced by \cite{donohoj1}. It can be
expressed in the form: $\lambda^*=(\lambda^*_j)_j$ where
$\lambda^*_j=c \sqrt{(\log n)/n}$ where $c>0$ denotes a large enough
constant. In the literature, several technics have been proposed to
determine the 'best' adaptive threshold.
There are, for instance, the RiskShrink and SureShrink methods (see
\cite{donoho4,donohoj1}), the cross-validation methods (see
\cite{nason}, \cite{weyrich} and \cite{jansen}), the methods based
on hypothesis tests (see \cite{abramovich} and \cite{abdj:06}), the
Lepski methods (see \cite{juditsky}) and the Bayesian methods (see
\cite{chipman} and \cite{abramovich2}). Most of them are described
in detailed in \cite{nason} and \cite{antoniadis}.

In the present paper, we propose to study the performances of an
adaptive wavelet estimator based on a convex combination of $\hat
f_{\lambda}$'s.
In the framework of nonparametric density estimation and bounded
regression estimation with random uniform design, we prove that, in
some sense, it is at least as good as the term-by-term thresholded
estimator $\hat f_{\lambda}$ defined with the 'best' threshold
$\lambda$. In particular, we show that this estimator is optimal, in
the minimax sense, over all Besov balls under the $L^2$ risk. The
proof is based on a non-adaptive minimax result proved by
\cite{delyon} and some powerful oracle inequality satisfied by
aggregation methods. There are two steps in our approach. A first
step, called the training step, where non-adaptive thresholded
wavelet estimators are constructed for different thresholds. A
second step, called learning step, where an aggregation scheme is
worked out to realize the adaptation to the smoothness.

The exact oracle inequality of Section 2 is given in a general
framework. Two aggregation procedures satisfy this oracle
inequality. The well known ERM (for Empirical Risk Minimization)
procedure (cf. \cite{vbook:98}, \cite{k:05} and references therein)
and an exponential weighting aggregation scheme, which has been
studied, among others, by \cite{bl:04}, \cite{bn:05}, \cite{lec:05},
\cite{lec4:06} and \cite{lec5:05}. There is a recursive version of
this scheme studied by \cite{catbook:01}, \cite{yang:00},
\cite{jntv:05} and \cite{jrt:06}. In the sequential prediction
problem, weighted average predictions with exponential weights have
been widely studied (cf. e.g. \cite{v:90} and \cite{cblbook:06}). A
recent result of \cite{lec6:06} shows that the ERM procedure is
suboptimal for strictly convex losses (which is the case for density
and regression estimation when the integrated squared risk is used).
Thus, in our case it is better to combine the $\hat{f}_\lambda$'s,
for $\lambda$ lying in a grid, using the aggregation procedure with
exponential weights than using the ERM procedure. Moreover, from a
computation point of view the aggregation scheme with exponential
weights does not require any minimization step contrarily to the ERM
procedure.

The paper is organized as follows. Section 2 presents general oracle
inequalities satisfied by two aggregation methods. Section 3
describes the main procedure of the study and investigates its
minimax performances over Besov balls for the $L^2$ risk.
All the proofs are postponed in the last section.

\section{Oracle Inequalities}
\subsection{Framework}\label{SubSectionFramework}
Let $(\cZ,\cT)$ a measurable space. Denote by $\cP$ the set of all
probability measures on $(\cZ,\cT)$. Let $F$ be a function from
$\cP$ with values in an algebra $\cF$. Let $Z$ be a random variable
with values in $\cZ$ and denote by $\pi$ its probability measure.
Let $D_n$ be a family of $n$ i.i.d. observations $Z_1,\ldots,Z_n$
having the common probability measure $\pi$. The probability measure
$\pi$ is unknown. Our aim is to estimate $F(\pi)$ from the
observations $D_n$.

In our estimation problem, we assume that we have access to an
"empirical risk". It means that there exists
$Q:\cZ\times\cF\longmapsto\mathbb{R}$ such that the risk of an
estimate $f\in\cF$ of $F(\pi)$ is of the form
$$A(f)=\mathbb{E}\left[ Q(Z,f)\right].$$ In what follows, we present
several statistical problems which can be written in this way. If
the minimum over all $f$ in $\cF$
$$A^*\egal\min_{f\in\cF}A(f)$$ is achieved by at least one function, we
denote by $f^*$ a minimizer in $\cF$. In this paper we will assume
that $\min_{f\in\cF}A(f)$ is achievable, otherwise we replace
$f^*$ by $f^*_n$, an element in $\cF$ satisfying $A(f^*_n)\leq
\inf_{f\in\cF}A(f)+n^{-1}.$

In most of the cases $f^*$ will be equal to our aim $F(\pi)$ up to
some known additive terms. We don't know the risk $A$, since $\pi$
is not available from the statistician, thus, instead of
minimizing $A$ over $\cF$ we consider an empirical version of $A$
constructed from the observations $D_n$. The main interest of such
a framework is that we have access to an empirical version of
$A(f)$ for any $f\in\cF$. It is denoted by
\begin{equation}\label{EmpiricalRisk}A_n(f)=\frac{1}{n}\sum_{i=1}^nQ(Z_i,f).\end{equation}

We exhibit three statistical models having the previous form of
estimation.

{\bf{Bounded Regression:}} Take $\cZ=\cX\times[0,1]$, where
$(\cX,\cA)$ is a measurable space, $Z=(X,Y)$ a couple of random
variables on $\cZ$, with probability distribution $\pi$,  such
that $X$ takes its values in $\cX$ and $Y$ takes its values in
$[0,1]$. We assume that the conditional expectation
$\mathbb{E}[Y|X]$ exists. In the regression framework, we want to
estimate the regression function
$$f^*(x)=\mathbb{E}\left[Y|X=x \right],\ \forall x\in\cX.$$
Usually, the variable $Y$ is not an exact function of $X$. Given
is an input $X\in\cX$, we are not able to predict the exact value
of the output $Y\in[0,1]$. This issue can be seen in the
regression framework as a noised estimation. It means that in each
spot $X$ of the input set, the predicted label $Y$ is concentrated
around $\mathbb{E}\left[Y|X \right]$ up to an additional noise
with null mean denoted by $\zeta$. The regression model can then
be written as
$$Y=\mathbb{E}\left[Y|X \right]+\zeta.$$
Take $\cF$ the set of all measurable functions from $\cX$ to
$[0,1]$. Define $||f||_{L^2(P^X)}^2=\int_{\cX}f^2(x)dP^X(x)$ for
all functions $f$ in $L^2(\cX,\cA,P^X)$ where $P^X$ is the
probability measure of $X$. Consider
\begin{equation}\label{LossRegression}Q((x,y),f)=(y-f(x))^2,\end{equation} for any
$(x,y)\in\cX\times\mathbb{R}$ and $f\in\cF$. Pythagore's Theorem
yields
$$A(f)=\mathbb{E}\left[Q((X,Y),f)\right]=||f^*-f||_{L^2(P^X)}^2+\mathbb{E}\left[\zeta^2 \right].$$
Thus $f^*$ is a minimizer of $A(f)$ and $A^*=\mathbb{E}[\zeta^2].$

{\bf{Density estimation:}} Let $(\cZ,\cT,\mu)$ be a measured
space. Let $Z$ be a random variable with values in $\cZ$ and
denote by $\pi$ its probability distribution. We assume that $\pi$
is absolutely continuous w.r.t. to $\mu$ and denote by $f^*$ one
version of the density. Consider $\cF$ the set of all density
functions on $(\cZ,\cT,\mu)$. We consider
$$Q(z,f)=-\log f(z),$$ for any $z\in\cZ$ and $f\in\cF$.
We have
$$A(f)=\mathbb{E}\left[Q(Z,f) \right]=K(f^*|f)-\int_{\cZ}\log
(f^*(z))d\pi(z).$$ Thus, $f^*$ is a minimizer of $A(f)$ and
$A^*=-\int_{\cZ}\log (f^*(z))d\pi(z)$.

Instead of using the Kullback-Leiber loss, one can use the
quadratic loss. For this setup, consider $\cF$ the set
$L^2(\cZ,\cT,\mu)$ of all measurable functions with an integrated
square.
Define\begin{equation}\label{LossDensity}Q(z,f)=\int_{\cZ}f^2d\mu-2f(z),\end{equation}for
any $z\in\cZ$ and $f\in\cF$. We have, for any $f\in\cF$,
$$A(f)=\mathbb{E}\left[Q(Z,f) \right]=||f^*-f||^{2}_{L^2(\mu)}-\int_{\cZ}
(f^*(z))^2d\mu(z).$$ Thus, $f^*$ is a minimizer of $A(f)$ and
$A^*=-\int_{\cZ}(f^*(z))^2d\mu(z)$.

{\bf{Classification framework:}} Let $(\cX,\cA)$ be a measurable
space. We assume that the space $\cZ=\cX\times\{-1,1\}$ is endowed
with an unknown probability measure $\pi$. We consider a random
variable $Z=(X,Y)$ with values in $\cZ$ with probability
distribution $\pi$. We denote by $P^X$ the marginal of $\pi$ on
$\cX$ and $\eta(x)=\mathbb{P}(Y=1|X=x)$ the conditional
probability function of $Y=1$ knowing that $X=x$. Denote by $\cF$
the set of all measurable functions from $\cX$ to $\mathbb{R}$.
Let $\phi$ be a function from $\mathbb{R}$ to $\mathbb{R}$. For
any $f\in\cF$ consider the $\phi-$risk
$$A(f)=\mathbb{E}[Q((X,Y),f)],$$  where the loss is given by
$Q((x,y),f)=\phi(yf(x))$for any $(x,y)\in\cX\times\{-1,1\}$.

Most of the time a minimizer $f^*$ of the $\phi-$risk $A$ over
$\cF$ or its sign is equal to the Bayes rule $f^*(x)={\rm
Sign}(2\eta(x)-1),\forall x\in\cX$ (cf. \cite{z:04}).

In this paper we obtain an oracle inequality in the general
framework described at the beginning of this Subsection. Then, we
use it in the density estimation and the bounded regression
frameworks. For applications of this oracle inequality in the
classification setup, we refer to \cite{lec4:06} and
\cite{lec:05}.

Now, we introduce an assumption which improve the quality of
estimation in our framework. This assumption has been first
introduced by \cite{mt:99}, for the problem of discriminant
analysis, and \cite{tsy:04}, for the classification problem. With
this assumption, parametric rates of convergence can be achieved,
for instance, in the classification problem (cf. \cite{tsy:04},
\cite{ss:04}).

{\bf Margin Assumption(MA):} {\it The probability measure $\pi$
satisfies the margin assumption MA($\kappa,c,\cF_0$), where
$\kappa\geq1,c>0$ and $\cF_0$ is a subset of $\cF$ if
$$\mathbb{E}[(Q(Z,f)-Q(Z,f^*))^2]\leq c(A(f)-A^*)^{1/\kappa}, $$
for any function $f\in\cF_0$.}

In the bounded regression setup, it is easy to see that any
probability distribution $\pi$ on $\cX\times[0,1]$ naturally
satisfies the margin assumption MA($1,16,\cF_1$), where $\cF_1$ is
the set of all measurable functions from $\cX$ to $[0,1]$. In
density estimation with the integrated squared risk, all
probability measures $\pi$ on $(\cZ,\cT)$ absolutely continuous
w.r.t. the measure $\mu$ with one version of its density a.s.
bounded by a constant $B\geq1$, satisfies the margin assumption
MA($1,16B^2,\cF_B$) where $\cF_B$ is the set of all non-negative
function $f\in L^2(\cZ,\cT,\mu)$ bounded by $B$.

Actually, the margin assumption is linked to the convexity of the
underlying loss. In density and regression estimation it is
naturally satisfied with the better margin parameter $\kappa=1$,
but, for non-convex loss (for instance in classification) this
assumption does not hold naturally (cf. \cite{lec6:06} for a
discussion on the margin assumption and for examples of losses
which does not satisfied naturally the margin assumption with
parameter $\kappa=1$).

\par

\subsection{Aggregation Procedures}
Let's work with the notations introduced in the beginning of the
previous Subsection. The aggregation framework considered, among
others, by \cite{jn:00}, \cite{yang:00},
\cite{catbook:01},\cite{n:00}, \cite{tsy:03}, \cite{bl:04},
\cite{b:04} is the following:  take $\cF_0$ a finite subset of
$\cF$, our aim is to mimic (up to an additive residual) the best
function in $\cF_0$ w.r.t. the risk $A$. For this, we consider two
aggregation procedures.

The Aggregation with Exponential Weights aggregate ({\bf AEW})
over $\cF_0$ is defined by
\begin{equation}\label{AEW}\tilde{f}_n^{(AEW)}\egal\sum_{f\in\cF_0}w^{(n)}(f)f,\end{equation} where
the exponential
  weights $w^{(n)}(f)$ are defined by
  \begin{equation}\label{weightsAEW}w^{(n)}(f)=
  \frac{\exp\left( -nA_n(f)\right)}{\sum_{g\in\cF_0}\exp\left(
 -nA_n(g)\right)},\quad \forall f\in\cF_0.\end{equation}

We consider the Empirical Risk Minimization procedure {\bf{(ERM)}}
over $\cF_0$ defined by
\begin{equation}\label{ERM}\tilde{f}_n^{(ERM)}\in {\rm Arg
}\min_{f\in\cF_0}A_n(f).\end{equation}

\par

\subsection{Oracle Inequalities}
In this Subsection we state an exact oracle inequality satisfied
by the ERM procedure and the AEW procedure (in the convex case) in
the general framework of the beginning of Subsection
\ref{SubSectionFramework}. From this exact oracle inequality we
deduce two others oracle inequalities in the density estimation
and the bounded regression framework. We introduce a quantity
which is going to be our residual term in the exact oracle
inequality. We consider
$$\gamma(n,M,\kappa,\cF_0,\pi,Q)=\left\{
    \begin{array}{ll}
        \left(\frac{\cB(\cF_0,\pi,Q)^{\frac{1}{\kappa}}\log M}{\beta_1n}\right)^{1/2} & \mbox{if }
            \cB(\cF_0,\pi,Q)\geq\left(\frac{\log
            M}{\beta_1 n}\right)^{\frac{\kappa}{2\kappa-1}}\\
        \left(\frac{\log M}{\beta_2n}\right)^{\frac{\kappa}{2\kappa-1}}&\mbox{otherwise, }\\
    \end{array}\right.$$where $\cB(\cF_0,\pi,Q)$ denotes
$\min_{f\in\cF_0}\left(A(f)-A^{*} \right)$, $\kappa\geq1$ is the
margin parameter, $\pi$ is the underlying probability measure, $Q$
is the loss function,
\begin{equation}\label{equaBeta1}\beta_1=\min\Big(\frac{\log2}{96cK},
\frac{3\sqrt{\log2}}{16K\sqrt{2}},\frac{1}{8(4c+K/3)},\frac{1}{576c}\Big).\end{equation}
and
\begin{equation}\label{Beta2}\beta_2=\min
\Big(\frac{1}{8},\frac{3\log2}{32K},\frac{1}{2(16c+K/3)},\frac{\beta_1}{2}\Big),\end{equation}
where the constant $c>0$ appears in MA($\kappa,c,\cF_0$).
\begin{thm}\label{TheoExactOracleIneq}
Consider the general framework introduced in the beginning of
Subsection \ref{SubSectionFramework}. Let $\cF_0$ denote a finite
subset of $M$ elements $f_1,\ldots,f_M$ in $\cF$, where $M\geq2$
is an integer. Assume that the underlying probability measure
$\pi$ satisfies the margin assumption MA($\kappa,c,\cF_0$) for
some $\kappa\geq1,c>0$ and $|Q(Z,f)-Q(Z,f^*)|\leq K$ a.s., for any
$f\in\cF_0$, where $K\geq1$ is a constant. The Empirical Risk
Minimization procedure (\ref{ERM}) satisfies
\begin{eqnarray*}
 \mathbb{E}[A(\tilde{f}_n^{(ERM)})-A^*]\leq
 \min_{j=1,\ldots,M}(A(f_j)-A^*)+4\gamma(n,M,\kappa,\cF_0,\pi,Q).
\end{eqnarray*}

 Moreover, if $f\longmapsto Q(z,f)$ is convex
for $\pi$-almost $z\in\cZ$, then the AEW procedure satisfies the
same oracle inequality as the ERM procedure.
\end{thm}
Now, we give two corollaries of Theorem \ref{TheoExactOracleIneq}
in the density estimation and bounded regression framework.
\begin{cor}\label{oraclerégression}
Consider the bounded regression setup. Let $f_1,\ldots,f_M$ be $M$
functions on $\cX$ with values in $[0,1]$. Let $\tilde{f}_n$
denote either the ERM or the AEW procedure. For $\beta_2$ defined
in (\ref{Beta2}) and for any $\epsilon>0$, we have
\begin{eqnarray*}
 \mathbb{E}[||f^*-\tilde{f}_n||_{L^2(P^X)}^2]  \leq  (1+\epsilon)
\min_{j=1,\ldots,M}(||f^*-f_j||_{L^2(P^X)}^2) +\frac{4\log
M}{\epsilon\beta_2 n}.
\end{eqnarray*}
\end{cor}
\begin{cor}\label{oracledensité}
Consider the density estimation framework. Assume that the
underlying density function $f^*$ to estimate is bounded by
$B\geq1$. Let $f_1,\ldots,f_M$ be $M$ functions bounded from above
and below by $B$. Let $\tilde{f}_n$ denote either the ERM or the AEW
procedure. For $\beta_2$ defined in (\ref{Beta2}) and any
$\epsilon>0$, we have\begin{equation}\label{EquaOracleDensity}
\mathbb{E}[||f^*-\tilde{f}_n||^{2}_{L^2(\mu)}]\leq (1+\epsilon)
\min_{j=1,\ldots,M}(||f^*-f_j||^{2}_{L^2(\mu)})+\frac{4\log
M}{\epsilon\beta_2n}.\end{equation}
\end{cor}
In both of the last Corollaries, the ERM and the AEW procedures can
both be used to mimic the best $f_j$ among the $f_j$'s.
Nevertheless, from a computational point of view the AEW procedure
does not require any minimization step contrarily to the ERM
procedure. Moreover, from a theoretical point of view the ERM
procedure can not mimic the best $f_j$ among the $f_j$'s as fast as
the cumulative aggregate with exponential weights (it is an average
of AEW procedures). For a comparison between these procedures we
refer to \cite{lec6:06}. The constants of aggregation multiplying
the residual term in Theorem \ref{TheoExactOracleIneq} and in both
of the following Corollaries come from the proof and are certainly
not optimal. We did not make any serious attempt to optimize them.


\par

\section{Multi-thresholding wavelet estimator}
In the present section, we propose an adaptive estimator
constructed from aggregation technics and wavelet thresholding
methods. For the density model and the regression model with
uniform random design, we show that it is optimal in the minimax
sense over a wide range of function spaces.
\subsection{Wavelets and Besov balls}
We consider an orthonormal wavelet basis generated by dilation and
translation of  a compactly supported "father" wavelet $\phi$ and
a compactly supported "mother" wavelet $\psi$. For the purposes of this paper, we use the
periodized wavelets bases on the unit interval. Let
$$\phi_{j,k}=2^{j/2}\phi (2^j x-k), \ \ \ \ \ \ \ \ \ \ \ \ \
\psi_{j,k}=2^{j/2}\psi(2^j x-k)$$be the elements of the wavelet
basis and
$$\phi^{per}_{j,k}(x)=\sum_{l\in \mathbb{Z}}\phi_{j,k}(x-l), \ \ \ \ \ \ \
 \ \ \ \psi^{per}_{j,k}(x)=\sum_{l\in
 \mathbb{Z}}\psi_{j,k}(x-l),$$there periodized versions, defined
 for any $x\in[0,1]$, $j\in\mathbb{N}$ and
 $k\in\{0,\ldots,2^j-1\}$.
There exists an integer $\tau$ such that the collection $\zeta$
defined by $\zeta = \{\phi^{per}_{j,k},  k =0,...,2^{\tau}-1; \
\psi^{per}_{j,k}, \ \ \ j=\tau,..., \infty, \ k=0,...,2^j-1\}$
constitutes an orthonormal basis of $L^2(\lbrack 0,1 \rbrack)$. In
what follows, the superscript "$per$" will be suppressed from the
notations for convenience. For any integer $l\geq\tau$, a
square-integrable function $f^*$ on $\lbrack 0,1 \rbrack$ can be
expanded into a wavelet series
$$f^*(x)=\sum_{k=0}^{2^{l}-1}\alpha_{l,k}\phi_{l,k}(x)+\sum_{j=l}^{\infty}\sum_{k=0}^{2^j-1}\beta_{j,k}\psi_{j,k}(x),$$
where $\alpha_{j,k}=\int_{0}^{1}f^*(x)\phi_{j,k}(x)dx$ and
$\beta_{j,k}=\int_{0}^{1}f^*(x)\psi_{j,k}(x)dx$. Further details
on wavelet theory can be found
in \cite{meyer} and \cite{daubechies}. 

Now, let us define the main function spaces of the study. Let
$M\in ( 0, \infty)$, $s\in (0, N )$, $ p \in \lbrack 1,\infty)$
and $q \in \lbrack 1,\infty)$. Let us set
$\beta_{\tau-1,k}=\alpha_{\tau,k}$. We say that a function $f^*$
belongs to the Besov balls $ B^s_{p,q}(M)$ if and only if the
associated wavelet coefficients satisfy
$$\Big[\sum_{j=\tau-1}^{\infty} \Big[2^{j(s+1/2-1/p )}\Big(\sum_{k=0}^{2^{j}-1}
|\beta_{j,k}|^p\Big)^{1/p}\Big]^{q}\Big]^{1/q}\le M, \ \ \ \ \  if
\ \ q\in \lbrack 1,\infty),
$$
with the usual modification if $q=\infty$. We work with the Besov
balls because of their exceptional expressive power. For a
particular choice of parameters $s$, $p$ and $q$, they contain the
H{\"o}lder and Sobolev balls (see \cite{meyer}).

\par

\subsection{Term-by-term thresholded estimator}\label{rob}
In this Subsection, we consider the estimation of an unknown
function $f^*$ in $L^2(\lbrack 0,1 \rbrack)$ from a general
situation. We only assume to have $n$ observations gathered in the
data set $D_n$ from which we are able to estimate the wavelet
coefficients $\alpha_{j,k}$ and $\beta_{j,k}$ of $f^*$ in the
basis $\zeta$. We denote by $\hat \alpha_{j,k}$ and $\hat
\beta_{j,k}$ such estimates. Finally, let us mention that all the
constants of our study are independent of $f^*$ and $n$.
\begin{defn}[Term-by-term thresholded estimator]
Let $j_1$ be an integer satisfying $(n/\log n)\le 2^{j_1}<2 (n/ \log n)$. For
any integer $l \geq\tau$, let
$\lambda=(\lambda_{l},...\lambda_{j_1})$ be a vector of positive
integers. Let us consider the estimator $\hat{f}_{\lambda}:
D_n\times \lbrack 0,1 \rbrack \rightarrow \mathbb{R}$ defined by
\begin{eqnarray} \label{hinz}
\hat{f}_{\lambda} (D_n,x)= \sum_{k=0}^{2^\tau-1}\hat \alpha_{\tau,k}
\phi_{\tau,k}(x)+\sum_{j=\tau}^{j_1}
\sum_{k=0}^{2^j-1}\Upsilon_{\lambda_j}(\hat \beta_{j,k})
\psi_{j,k}(x),
\end{eqnarray}where for all
$u\in(0,\infty)$ the operator $\Upsilon_u$ is such that there
exist two constants $C_1,C_2>0$ satisfying
\begin{eqnarray}\label{ongle}
|\Upsilon_{u}(x ) -y |^2\le C_1 ( \min(y,C_2 u)^2+(|x -y|^2) \1_{
\left \lbrace |x -y|\ge 2^{-1}u\right \rbrace }),
\end{eqnarray}for any $x\in\mathbb{R} $ and $y\in\mathbb{R}$.
\end{defn}
The inequality (\ref{ongle}) holds for the hard thresholding rule
$\Upsilon_{u}^{hard}(x ) = x\1_{\left \lbrace |x|\geqslant u \right
\rbrace}$, the soft thresholding rule
$\Upsilon_{u}^{soft}(x)=sign(x)(|x|-u)\1_{\left \lbrace |x
|\geqslant u \right \rbrace}$ (see \cite{donohoj1}, \cite{donoho001}
and \cite{delyon}) and the non-negative garrote thresholding rule
$\Upsilon_{u}^{NG}(x) =\left(x -{u^2}/{x}\right) \1_{\left \lbrace
|x |\geqslant u \right \rbrace}$ (see \cite{gao}).

If we consider the minimax point of view over Besov balls under the
integrated squared risk, then \cite{delyon} makes the conditions on
$\hat \alpha_{j,k}$, $\hat \beta_{j,k}$ and the threshold $\lambda$
such that the estimator $\hat f_{\lambda}(D_n,.)$ defined by
(\ref{hinz}) is optimal for numerous statistical models. This result
is recalled in Theorem \ref{guy} below.
\begin{thm}[Delyon and Juditsky (1996)]\label{guy}
Let us consider the general statistical framework described in the
beginning of the present section. Suppose that the two following
assumptions hold.
\begin{itemize}
\item \textit{Moments inequality:} There exists a constant $C>0$
such that, for any $j\in \{\tau-1,...,j_1\}$, $k\in \{0,...,2^{j}-1\}$
and $n$ large enough, we have
 \begin{eqnarray}\label{moment}
  \mathbb{E} (|\hat \beta_{j,k}-\beta_{j,k}|^{4}) \leq
  Cn^{-2},\mbox{ where we take
  }\hat \beta_{\tau-1,k}=\hat \alpha_{\tau,k}.
\end{eqnarray}
\item \textit{Large deviation inequality:} There exist two
constants $C>0$ and $\rho_*>0$ such that, for any $a, j\in
\{\tau,...,j_1\}$, $k\in \{0,...,2^{j}-1\} $ and $n$ large enough,
we have
\begin{eqnarray}\label{deviation}
 \mathbb{P}\left(2\sqrt{n}|\hat
\beta_{j,k}-\beta_{j,k}|\geq \rho_* \sqrt{a}\right) \leqslant C
2^{-4a}.
\end{eqnarray}
\end{itemize}
Let us consider the term-by-term thresholded estimator $\hat
f_{v_{j_s}} (D_n,.) $ defined by (\ref{hinz}) with the threshold
$$v_{j_s}=( \rho_*(j-j_s)_{+} ) _{j=\tau,...,j_1},$$ where $j_s$ is
an integer such that $n^{1/(1+2s)}\le 2^{j_s}< 2 n^{1/(1+2s)} $.
Then, there exists a constant $C>0$ such that, for any $p\in
\lbrack 1,\infty\rbrack$, $s\in (1/p, N \rbrack$, $q\in \lbrack
1,\infty\rbrack$ and $n$ large enough, we have:
\begin{eqnarray*}
\sup_{f\in {B}^s_{p,q}(L)}  \mathbb{E}[\Vert \hat f_{v_{j_s}}
(D_n,.) -f^*\Vert^{2}_{L^2(\lbrack 0,1\rbrack)} ]\leqslant C n^{-2s/(2s+1)}.
\end{eqnarray*}
\end{thm}
The rate of convergence $V_n=n^{-2s/(1+2s)}$ is minimax for
numerous statistical models, where $s$ is a regularity parameter.
For the density model and the regression model with uniform design,
we refer the reader to \cite{delyon} for further details about the
choice of the estimator $\hat \beta_{j,k}$ and the value of the thresholding constant $\rho_*$.
Starting from this non-adaptive result, we use aggregation methods
to construct an adaptive estimator at least at good in the minimax sense as $\hat f_{v_{j_s}} (D_n,.)$.

\subsection{Multi-thresholding estimator}
Let us divide our observations $D_{n}$ into two disjoint subsamples
$D_m$, of size $m$, made of the first $m$ observations and
$D^{(l)}$, of size $l$, made of the last remaining observations,
where we take$$l=\left\lceil {n}/{\log n}\right\rceil \mbox{ and }
m=n-l.$$ The first subsample $D_m$, sometimes called "training
sample", is used to construct a family of estimators (in our case
this is thresholded estimators) and the second subsample $D^{(l)}$,
called the "training sample", is used to construct the weights of
the aggregation procedure.
\begin{rem} From a theoretical point of view we can take $m=l$ which
means that we use as many observations for the estimation step as
for the learning step. But, in practice it is better to use a
greater part of the observations for the construction of the
estimators and the last observations for the aggregation
procedure, because if the basis estimators that we aggregate, are
not good, then the obtained aggregate is likely to be as bad as
the prior estimators. Another interesting thing is that we can
split the whole sample $D_n$ in many different ways. For instance
we can take $m$ observations randomly in $D_n$ to form the
training subsample and the last remaining observations for the
learning subsample. We can also take an average of different
aggregates constructed from different splits of the initial sample
$D_n$ and by a simple argument of convexity it is easy to prove
that the averaged aggregate has a better risk than the others
aggregates constructed only from one split.
\end{rem}
\begin{defn} Let us consider the term-by-term
thresholded estimator described in (\ref{hinz}). Assume that we want
to estimate a function $f^*$ from $[0,1]$ with values in $[a,b]$.
Consider the projection function
\begin{equation}\label{EquaProjFunc}
h_{a,b}    (y)=\max(a,\min(y,b)), \forall y\in\mathbb{R}.
\end{equation} We
define the {\bf{multi-thresholding estimator}} $\tilde{f}_n: [0,1]
\rightarrow [a,b]$ at a point $x\in[0,1]$ by the following
aggregate
\begin{eqnarray}\label{multi-thresholding}
\tilde f_n(x) =\sum_{u \in \Lambda_n} w^{(l)}(h_{a,b}   (\hat
f_{v_u}(D_m,.) ))
 h_{a,b}   (\hat f_{v_u}(D_m,x)),
\end{eqnarray}
where $\Lambda_n=\{0,...,\log n\}$,
$v_u=(\rho(j-u)_{+})_{j=\tau,...,j_1},\forall u \in \Lambda_n$ and
$\rho$ is a positive constant depending on the model worked out and
 $$w^{(l)}(h_{a,b}   (\hat f_{v_u}(D_m,.)))=
\frac{\exp\left( -lA^{(l)}(h_{a,b}  (\hat
f_{v_u}(D_m,.)))\right)}{\sum_{\gamma \in \Lambda_n}\exp\left(
 -lA^{(l)}(h_{a,b}  (\hat
f_{v_{\gamma}}(D_m,.)))\right)},\quad \forall u \in \Lambda_n,$$
where $A^{(l)}(f)=\frac{1}{l}\sum_{i=m+1}^n Q(Z_i,f)$ is the
empirical risk constructed from the $l$ last observations, for any
function $f$ and for the choice of a loss function $Q$ depending
on the model considered (cf. (\ref{LossRegression}) and
(\ref{LossDensity}) for examples).
\end{defn}

The principle of the construction of the multi-thresholding
estimator $\tilde f_n$ is to use aggregation technics to easily
construct an adaptive optimal estimator of $f^*$. It realizes a
kind of 'adaptation to the threshold' by selecting the best
threshold $v_u$ for $u$ describing the set $\Lambda_n$. Since we
know that there exists an element in $\Lambda_n$ depending on the
regularity of $f^*$ such that the non-adaptive estimator $\hat
f_{v_u}(D_m,.)$ is optimal in the minimax sense (see Theorem
\ref{guy}), the multi-thresholding estimator is optimal
independently of the regularity of $f^*$.

\par

\section{Performances of the multi-thresholding estimator}
This section is devoted to the minimax performances of the
multi-thresholding estimator defined in (\ref{multi-thresholding})
under the $L^2([0,1])$ risk over Besov balls. Firstly, we consider
the framework of the density model. Secondly, we focus our attention
on the bounded regression with uniform random design. Finally, we
compare these results with some well-known wavelet thresholded
procedures.

\subsection{Density model}
In the density estimation model, Theorem \ref{multiperf} below
investigates rates of convergence achieved by the
multi-thresholding estimator (defined by
(\ref{multi-thresholding})) under the $L^2([0,1])$ risk over Besov
balls.
\begin{thm}\label{multiperf}
Let us consider the problem of estimating $f^*$ from the density
model. Assume that there exists $B\geq1$ such that the underlying
density function $f^*$ to estimate is bounded by $B$. Let us
consider the multi-thresholding estimator defined in
(\ref{multi-thresholding}) where we take $a=0,b=B$, $\rho$ such
that
$$\frac{\rho^2}{8B+(8\rho/(3\sqrt{2}))(\Vert \psi\Vert_{\infty}+B)}\ge 4(\log 2)$$
and
\begin{eqnarray}\label{est1}
\hat \alpha_{j,k}=\frac{1}{n}\sum_{i=1}^n \phi_{j,k}(X_i), \ \ \ \
\ \ \ \ \ \ \ \ \ \hat \beta_{j,k}=\frac{1}{n}\sum_{i=1}^n
\psi_{j,k}(X_i).
\end{eqnarray}
Then, there exists a constant $C>0$ such that
$$ \sup_{f^*\in {{B}}^{s} _{p,q}  (L)}  \mathbb{E}[\Vert
\tilde{f}_n -f^*\Vert^{2}_{L^2([0,1])}] \leqslant
Cn^{-2s/(2s+1)},$$ for any $p\in [1,\infty]$, $s\in(p^{-1},N]$,
$r\in [1,\infty]$ and integer $n$.
\end{thm}
The rate of convergence $V_n=n^{-2s/(1+2s)}$ is minimax over
${{B}}^{s}_{p,q} (L)$. Further details about the minimax rate of
convergence over Besov balls under the $L^2([0,1])$ risk for the
density model can be found in \cite{delyon} and \cite{hardle}. For
further details about the density estimation via adaptive wavelet
thresholded estimators, see \cite{donoho2}, \cite{delyon} and
\cite{picard}. See also \cite{silverman15} for a practical study.

\par

\subsection{Bounded regression}
In the framework of the bounded regression model with uniform
random design, Theorem \ref{multiperf2} below investigates the
rate of convergence achieved by the multi-thresholding estimator
defined by (\ref{multi-thresholding}) under the $L^2([0,1])$ risk
over Besov balls.
\begin{thm}\label{multiperf2}
Let us consider the problem of estimating the regression function
$f^*$ in the bounded regression model with random uniform design.
Let us consider the multi-thresholding estimator
(\ref{multi-thresholding}) with $\rho$ such that
$$\frac{\rho^2}{8+(8\rho/(3\sqrt{2}))(\Vert \psi\Vert_{\infty}+1)}\ge 4(\log 2)$$
 and
\begin{eqnarray}\label{est2}
\hat \alpha_{j,k}=\frac{1}{n}\sum_{i=1}^n Y_i\phi_{j,k}(X_i), \ \
\ \ \ \ \ \ \ \ \ \ \ \hat \beta_{j,k}=\frac{1}{n}\sum_{i=1}^n Y_i
\psi_{j,k}(X_i).
\end{eqnarray}
 Then, there exists a constant $C>0$ such that, for any
$p\in [ 1,\infty]$, $s\in (p^{-1},N]$, $q\in \lbrack
1,\infty\rbrack$ and integer $n$, we have
$$ \sup_{f^*\in B^s_{p,q} (L)}  \mathbb{E}[\Vert \tilde{f}_n  -f^*\Vert_{L^2([0,1])}^2]
\leqslant Cn^{-2s/(2s+1)}.$$
\end{thm}

The rate of convergence $V_n=n^{-2s/(1+2s)}$ is minimax over
${{B}}^{s}_{p,q} (L)$. The multi-thresholding estimator has better
minimax properties than several other wavelet estimators developed
in the literature. To the authors's knowledge, the result obtained,
for instance, by the hard thresholded estimator (see
\cite{donohoj1}), by the global wavelet block thresholded estimator
(see \cite{kerk}), by the localized wavelet block thresholded
estimator (see \cite{caittt,cai4,cait1}, \cite{hall,hall2},
\cite{e:99,e:00}, \cite{chicken1} and \cite{chicken3}) and, in
particular, the penalized Blockwise Stein method (see
\cite{cavalier41}) are worse than the one obtained by the
multi-thresholding estimator and stated in Theorems \ref{multiperf}
and \ref{multiperf2}. This is because, on the difference of those
works, we obtain the optimal rate of convergence without any extra
logarithm factor.

In fact, the multi-thresholding estimator has similar minimax
performances than the empirical Bayes wavelet methods (see
\cite{zhang} and \cite{js}) and several term-by-term wavelet
thresholded estimators defined with a random threshold (see
\cite{juditsky} and \cite{birge}).

Finally, it is important to mention that the multi-thresholding
estimator does not need any minimization step and is relatively easy
to implement.

\section{Proofs}
{\bf{Proof of Theorem \ref{TheoExactOracleIneq}.}} We recall the
notations of the general framework introduced in the beginning of
Subsection \ref{SubSectionFramework}. Consider a loss function
$Q:\cZ\times\cF\longmapsto\mathbb{R}$, the risk
$A(f)=\mathbb{E}[Q(Z,f)]$, the minimum risk
$A^*=\min_{f\in\cF}A(f)$, where we assume, w.o.l.g, that it is
achieved by an element $f^*$ in $\cF$ and the empirical risk
$A_n(f)=(1/n)\sum_{i=1}^nQ(Z_i,f)$, for any $f\in\cF$. The
following proof is a generalization of the proof of Theorem 1 in
\cite{lec5:05}.

We first start by a 'linearization' of the risk. Consider the
convex set
$$\cC=\Big\{(\theta_1,\ldots,\theta_M):\theta_j\geq0 \mbox{ and } \sum_{j=1}^M
\theta_j=1\Big\}$$and define the following functions on $\cC$
$$\tilde{A}(\theta)\egal\sum_{j=1}^M \theta_j A(f_j) \mbox{ and } \tilde{A}_n(\theta)\egal\sum_{j=1}^M \theta_jA_n(f_j)$$
which are linear versions of the risk $A$ and its empirical
version $A_n$.

Using the Lagrange method of optimization we find that the
exponential weights $w\egal(w^{(n)}(f_j))_{1\leq j\leq M}$ are the
unique solution of the minimization problem
$$\min\Big(\tilde{A}_n(\theta)+\frac{1}{n}\sum_{j=1}^M \theta_j\log\theta_j:
 (\theta_1,\ldots,\theta_M)\in\cC \Big),$$
where we use the convention $0\log0=0$. Take
$\hat\jmath\in\{1,\ldots,M\}$ such that
$A_n(f_{\hat\jmath})=\min_{j=1,\ldots,M}A_n(f_j)$. The vector of
exponential weights $w$ satisfies
$$\tilde{A}_n(w)\leq \tilde{A}_n(e_{\hat\jmath})+\frac{\log M}{n},$$
where $e_j$ denotes the vector in $\cC$ with $1$ for $j$-th
coordinate (and $0$ elsewhere).

Let $\epsilon>0$. Denote by $\tilde{A}_\cC$ the minimum
$\min_{\theta\in\cC}\tilde{A}(\theta)$. We consider the subset of
$\cC$
$$\cD\egal\left\{\theta\in\cC:\tilde{A}(\theta)>\tilde{A}_\cC+2\epsilon\right\}.$$
Let $x>0$. If
$$\sup_{\theta\in\cD}\frac{\tilde{A}(\theta)-A^*-(\tilde{A}_n(\theta)-A_n(f^*))}{\tilde{A}(\theta)-A^*+x}\leq
\frac{\epsilon}{\tilde{A}_\cC-A^*+2\epsilon+x},$$then for any
$\theta\in\cD$, we have $$\tilde{A}_n(\theta)-A_n(f^*)\geq
\tilde{A}(\theta)-A^*-\frac{\epsilon(\tilde{A}(\theta)-A^*+x)}{(\tilde{A}_\cC-A^*+2\epsilon+x)}\geq
\tilde{A}_\cC-A^*+\epsilon,$$ because $\tilde{A}(\theta)-A^*\geq
\tilde{A}_\cC-A^*+2\epsilon$. Hence,
\begin{eqnarray}\label{equaPremIneg}\lefteqn{\mathbb{P}\left[\inf_{\theta\in\cD}
\left(\tilde{A}_n(\theta)-A_n(f^*)\right)<\tilde{A}_\cC-A^*+\epsilon
\right]}\nonumber\\ & \leq &
\mathbb{P}\left[\sup_{\theta\in\cD}\frac{\tilde{A}(\theta)-A^*-(\tilde{A}_n(\theta)-A_n(f^*))}
{\tilde{A}(\theta)-A^*+x}>
\frac{\epsilon}{\tilde{A}_\cC-A^*+2\epsilon+x}
\right].\end{eqnarray}

Observe that a linear function achieves its maximum over a convex
polygon at one of the vertices of the polygon. Thus, for
$j_0\in\{1,\ldots,M\}$ such that
$\tilde{A}(e_{j_0})=\min_{j=1,\ldots,M}\tilde{A}(e_j)\
(=\min_{j=1,\ldots,M} A(f_j))$, we have
$\tilde{A}(e_{j_0})=\min_{\theta\in\cC}\tilde{A}(\theta)$. We obtain
the last inequality by linearity of $\tilde{A}$ and the convexity of
$\cC$. Let $\hat{w}$ denotes either the exponential weights $w$ or
$e_{\hat\jmath}$. According to (\ref{equaPremIneg}), We have
$$\tilde{A}(\hat{w})\leq
\min_{j=1,\ldots,M}\tilde{A}_n(e_j)+\frac{\log M}{n} \leq
\tilde{A}_n(e_{j_0})+\frac{\log M}{n}$$ So, if
$\tilde{A}(\hat{w})>A_{\cC}+2\epsilon$ then $\hat{w}\in\cD$ and
thus, there exists $\theta\in\cD$ such that
$\tilde{A}_n(\theta)-\tilde{A}_n(f^*)\leq
\tilde{A}_n(e_{j_0})-\tilde{A}_n(f^*)+(\log M)/n$. Hence, we have
\begin{eqnarray*}
\lefteqn{\mathbb{P}\left[\tilde{A}(\hat{w})>\tilde{A}_{\cC}+2\epsilon
\right]\leq
\mathbb{P}\left[\inf_{\theta\in\cD}\tilde{A}_n(\theta)-A_n(f^*)
\leq \tilde{A}_n(e_{j_0})-A_n(f^*)+\frac{\log M}{n} \right]}\\
& \leq & \mathbb{P}\left[\inf_{\theta\in\cD}
\tilde{A}_n(\theta)-A_n(f^*)< \tilde{A}_\cC-A^*+\epsilon
\right]\\&&+\mathbb{P}\left[\tilde{A}_n(e_{j_0})-A_n(f^*)\geq
\tilde{A}_\cC-A^*+ \epsilon-\frac{\log M}{n} \right]\\
& \leq &
\mathbb{P}\left[\sup_{\theta\in\cC}\frac{\tilde{A}(\theta)-A^*-(\tilde{A}_n(f)-A_n(f^*))}
{\tilde{A}(\theta)-A^*+x}>
\frac{\epsilon}{\tilde{A}_\cC-A^*+2\epsilon+x} \right]\\
& & +\mathbb{P}\left[\tilde{A}_n(e_{j_0})-A_n(f^*)\geq
\tilde{A}_\cC-A^*+ \epsilon-\frac{\log M}{n} \right].
\end{eqnarray*}
If we assume that
$$\sup_{\theta\in\cC}\frac{\tilde{A}(\theta)-A^*-(\tilde{A}_n(\theta)-A_n(f^*))}{\tilde{A}(\theta)-A^*+x}>
\frac{\epsilon}{\tilde{A}_\cC-A^*+2\epsilon+x},$$ then, there
exists
$\theta^{(0)}=(\theta_1^{(0)},\ldots,\theta_M^{(0)})\in\cC$, such
that
$$\frac{\tilde{A}(\theta^{(0)})-A^*-(\tilde{A}_n(\theta^{(0)})-A_n(f^*))}{\tilde{A}(\theta^{(0)})-A^*+x}>
\frac{\epsilon}{\tilde{A}_\cC-A^*+2\epsilon+x}.$$ The linearity of
 $\tilde{A}$ yields
$$\frac{\tilde{A}(\theta^{(0)})-A^*-(\tilde{A}_n(\theta^{(0)})-A_n(f^*))}{\tilde{A}(\theta^{(0)})-A^*+x}=
\frac{\sum_{j=1}^M
\theta_j^{(0)}[A(f_j)-A^*-(A_n(f_j)-A_n(f^*))}{\sum_{j=1}^M\theta_j^{(0)}[A(f_j)-A^*+x]}$$
and since, for any  numbers $a_1,\ldots,a_M$ and positive numbers
$b_1,\ldots,b_M$, we have
$$\frac{\sum_{j=1}^Ma_j}{\sum_{j=1}^Mb_j}\leq
\max_{j=1,\ldots,M}\left(\frac{a_j}{b_j} \right),$$ then, we
obtain
$$\max_{j=1,\ldots,M}\frac{A(f_j)-A^*-(A_n(f_j)-A_n(f^*))}{A(f_j)-A^*+x}>
\frac{\epsilon}{A_{\cF_0}-A^*+2\epsilon+x},$$ where
$A_{\cF_0}\egal \min_{j=1,\ldots,M}A(f_j)\ (=\tilde{A}_\cC)$.

Now, we use the relative concentration inequality of Lemma
\ref{LemDevRela} to obtain
\begin{eqnarray*}
\lefteqn{\mathbb{P}\left[\max_{j=1,\ldots,M}\frac{A(f_j)-A^*-(A_n(f_j)-A_n(f^*))}{A(f_j)-A^*+x}>
\frac{\epsilon}{A_{\cF_0}-A^*+2\epsilon+x} \right]}\\
& \leq & M\left(
1+\frac{4c(A_{\cF_0}-A^*+2\epsilon+x)^2x^{1/\kappa}}{n(\epsilon
x)^2}\right)
\exp\left(-\frac{n(\epsilon x)^2}{4c(A_{\cF_0}-A^*+2\epsilon+x)^2x^{1/\kappa}} \right)\\
& & + M\left(1+\frac{4K(A_{\cF_0}-A^*+2\epsilon+x)}{3n\epsilon x}
\right)\exp\left(-\frac{3n\epsilon
x}{4K(A_{\cF_0}-A^*+2\epsilon+x)} \right).\end{eqnarray*} Using
the margin assumption MA($\kappa,c,\cF_0$) to upper bound the
variance term and applying Bernstein's inequality, we get
\begin{eqnarray*}\mathbb{P}\Big[A_n(f_{j_0})&-&A_n(f^*) \geq  A_{\cF_0}-A^*+
\epsilon-\frac{\log M}{n}\Big]\\
 &&\leq   \exp\left(-\frac{n(\epsilon-(\log M)
/n)^2}{2c(A_{\cF_0}-A^*)^{1/\kappa}+(2K/3)(\epsilon-(\log M)/n)}
\right),\end{eqnarray*} for any $\epsilon>(\log M)/n$. From now,
we take $x=A_{\cF_0}-A^*+2\epsilon$, then, for any $(\log M)/n <
\epsilon <1$, we have
\begin{eqnarray*}\lefteqn{\mathbb{P}\left(\tilde{A}(\hat{w})>A_{\cF_0}+2\epsilon
\right)\leq \exp\left(-\frac{n(\epsilon-\log M
/n)^2}{2c(A_{\cF_0}-A^*)^{1/\kappa}+(2K/3)(\epsilon-(\log M)/n)}
\right)}\\
&+&M\left(1+\frac{32c(A_{\cF_0}-A^*+2\epsilon)^{1/\kappa}}{n\epsilon^2}
\right)\exp\left(-
\frac{n\epsilon^2}{32c(A_{\cF_0}-A^*+2\epsilon)^{1/\kappa}}\right)\\
&+&M\left(1+\frac{32}{3n\epsilon}\right)\exp\left(-\frac{3n\epsilon}{32}
\right).\end{eqnarray*}If $\hat{w}$ denotes $e_{\hat\jmath}$ then,
$\tilde{A}(\hat{w})=\tilde{A}(e_{\hat\jmath})=A(\tilde{f}^{(ERM)})$.
If $\hat{w}$ denotes the vector of exponential weights $w$ and if
$f\longmapsto Q(z,f)$ is convex for $\pi$-almost $z\in\cZ$, then,
$\tilde{A}(\hat{w})=\tilde{A}(w)\geq A(\tilde{f}_n^{(AEW)})$. If
$f\longmapsto Q(z,f)$ is assumed to be convex for $\pi$-almost
$z\in\cZ$ then, let $\tf$ denote either the ERM procedure or the
AEW procedure, otherwise, let $\tf$ denote the ERM procedure
$\tilde{f}_n^{(ERM)}$. We have for any  $2(\log M)/n<u<1$,
\begin{equation}\label{equaInit}
\mathbb{E}[A(\tf)-A_{\cF_0}]\leq\mathbb{E}\left[\tilde{A}(\hat{w})-A_{\cF_0}
\right]\leq 2u+2
\int_{u/2}^1\left[T_1(\epsilon)+M(T_2(\epsilon)+T_3(\epsilon))
\right]d\epsilon,\end{equation} where
$$T_1(\epsilon)=\exp\left(-\frac{n(\epsilon-(\log M)
/n)^2}{2c(A_{\cF_0}-A^*)^{1/\kappa}+(2K/3)(\epsilon-(\log M)/n)}
\right),$$ $$
T_2(\epsilon)=\left(1+\frac{16c(A_{\cF_0}-A^*+2\epsilon)^{1/\kappa}}{n\epsilon^2}
\right)\exp\left(-
\frac{n\epsilon^2}{16c(A_{\cF_0}-A^*+2\epsilon)^{1/\kappa}}\right)$$
and
$$T_3(\epsilon)=\left(1+\frac{8K}{3n\epsilon}\right)\exp\left(-\frac{3n\epsilon}{8K}
\right).$$

We recall that  $\beta_1$ is defined in (\ref{equaBeta1}).
Consider separately the following cases ($C1$) and ($C2$).

\noindent{($C1$)\underline{ The case $A_{\cF_0}-A^*\geq ((\log
M)/(\beta_1 n))^{\kappa/(2\kappa-1)}$.}}

Denote by $\mu(M)$ the unique solution of $\mu_0=3M\exp(-\mu_0)$.
Then, clearly $(\log M)/2 \leq \mu(M)\leq \log M$. Take $u$ such
that
$$(n\beta_1u^2)/(A_{\cF_0}-A^*)^{1/\kappa}=\mu(M).$$ Using the
definition of case ($1$) and of $\mu(M)$ we get $u\leq
A_{\cF_0}-A^*$. Moreover, $u\geq 4\log M/n$, then
\begin{eqnarray*}
\int_{u/2}^{1}T_1(\epsilon)d\epsilon  &\leq&
\int_{u/2}^{(A_{\cF_0}-A^*)/2}\exp\left(-\frac{n(\epsilon/2)^2}{(2c+K/6)(A_{\cF_0}-A^*)^{1/\kappa}}\right)d\epsilon\\
&&+\int_{(A_{\cF_0}-A^*)/2}^{1}\exp\left(-\frac{n(\epsilon/2)^2}{(4c+K/3)\epsilon^{1/\kappa}}\right)d\epsilon.
\end{eqnarray*}
Using Lemma \ref{intalpha} and the inequality $u\leq
A_{\cF_0}-A^*$, we obtain
\begin{equation}\label{equa11}\int_{u/2}^{1}T_1(\epsilon)d\epsilon\leq
\frac{8(4c+K/3)(A_{\cF_0}-A^*)^{1/\kappa}}{nu}\exp\left(-
\frac{nu^2}{8(4c+K/3)(A_{\cF_0}-A^*)^{1/\kappa}}\right).
\end{equation}
We have $16c(A_{\cF_0}-A^*+2u)\leq nu^2$ thus, using Lemma
\ref{intalpha}, we get
\begin{eqnarray}\label{equa12}
\int_{u/2}^1T_2(\epsilon)d\epsilon &\leq&
2\int_{u/2}^{(A_{\cF_0}-A^*)/2}\exp\left(-\frac{n\epsilon^2}{64c(A_{\cF_0}-A^*)^{1/\kappa}}\right)d\epsilon
\nonumber\\
&&+2\int_{(A_{\cF_0}-A^*)/2}^{1}\exp\left(-\frac{n\epsilon^{2-1/\kappa}}{128c}
 \right)d\epsilon\nonumber\\
&\leq&\frac{2148c(A_{\cF_0}-A^*)^{1/\kappa}}{nu}\exp\left(-\frac{nu^2}{2148c(A_{\cF_0}-A^*)^{1/\kappa}}\right).
\end{eqnarray}
We have $16(3n)^{-1}\leq u\leq A_{\cF_0}-A^*$, thus,
\begin{equation}\label{equa13}\int_{u/2}^{1}T_3(\epsilon)d\epsilon
\leq\frac{16K(A_{\cF_0}-A^*)^{1/\kappa}}{3nu}\exp\left(-\frac{3nu^2}
{16K(A_{\cF_0}-A^*)^{1/\kappa}}\right).\end{equation} From
(\ref{equa11}), (\ref{equa12}), (\ref{equa13}) and
(\ref{equaInit}) we obtain
$$\mathbb{E}\left[A(\tf)-A_{\cF_0} \right]\leq
2u+6M\frac{(A_{\cF_0}-A^*)^{1/\kappa}}{n\beta_1u}
\exp\left(-\frac{n\beta_1u^2}{(A_{\cF_0}-A^*)^{1/\kappa}}
\right).$$ The definition of $u$ leads to
$\mathbb{E}\left[A(\tf)-A_{\cF_0} \right]\leq
4\sqrt{\frac{(A_{\cF_0}-A^*)^{1/\kappa}\log M}{n\beta_1}}.$

\noindent{($C2$)\underline{The case $A_{\cF_0}-A^*\leq ((\log
M)/(\beta_1 n))^{\kappa/(2\kappa-1)}$.}}

We now choose  $u$ such that $n\beta_2
u^{(2\kappa-1)/\kappa}=\mu(M)$, where $\mu(M)$ denotes the unique
solution of $\mu_0=3M\exp(-\mu_0)$ and $\beta_2$ is defined in
(\ref{Beta2}). Using the definition of case (2) and of $\mu(M)$ we
get $u\geq A_{\cF_0}-A^*$ (since $\beta_1\geq 2\beta_2$). Using
the fact that $u>4\log M/n$ and Lemma \ref{intalpha}, we have
\begin{equation}\label{equa21}\int_{u/2}^{1}T_1(\epsilon)d\epsilon\leq
\frac{2(16c+K/3)}{nu^{1-1/\kappa}}\exp\left(-\frac{3nu^{2-1/\kappa}}{2(16c+K/3)}
\right).\end{equation} We have $u\geq
(128c/n)^{\kappa/(2\kappa-1)}$ and using Lemma \ref{intalpha}, we
obtain
\begin{equation}\label{equa22}\int_{u/2}^{1}T_2(\epsilon)d\epsilon\leq
\frac{256c}{nu^{1-1/\kappa}}\exp\left(-\frac{nu^{2-1/\kappa}}{256c}
\right).\end{equation} Since $u>16K/(3n)$ we have
\begin{equation}\label{equa23}\int_{u/2}^{1}T_3(\epsilon)d\epsilon\leq
\frac{16K}{3nu^{1-1/\kappa}}\exp\left(-\frac{3nu^{2-1/\kappa}}{16K}
\right).\end{equation} From (\ref{equa21}), (\ref{equa22}),
(\ref{equa23}) and (\ref{equaInit}) we obtain
$$\mathbb{E}\left[A(\tf)-A_{\cF_0} \right]\leq
2u+6M\frac{\exp\left(-n\beta_2 u^{(2\kappa-1)/\kappa}
\right)}{n\beta_2 u^{1-1/\kappa}}.$$  The definition of $u$ yields
$\mathbb{E}\left[A(\tf)-A_{\cF_0} \right]\leq 4 \left(\frac{\log
M}{n\beta_2} \right)^{\frac{\kappa}{2\kappa-1}}.$ This completes
the proof.


\begin{lem}\label{LemDevRela} Consider the framework introduced in the
beginning of Subsection \ref{SubSectionFramework}. Let
$\cF_0=\{f_1,\ldots,f_M\}$ be a finite subset of $\cF$. We assume
that $\pi$ satisfies MA($\kappa,c,\cF_0$), for some $\kappa\geq 1,
c>0$ and $|Q(Z,f)-Q(Z,f^*)|\leq K$ a.s., for any $f\in\cF_0$,
where $K\geq1$ is a constant. We have for any positive numbers
$t,x$ and any integer $n$
\begin{eqnarray*} \lefteqn{\mathbb{P}\left[
\max_{f\in\cF}\frac{A(f)-A_n(f)-(A(f^*)-A_n(f^*))}{A(f)-A^*+x}>t\right]}\\
&  &\leq M\left(\left(1+\frac{4cx^{1/\kappa}}{n(tx)^2}
\right)\exp\left(
-\frac{n(tx)^2}{4cx^{1/\kappa}}\right)+\left(1+\frac{4K}{3ntx}
\right)\exp\left(-\frac{3ntx}{4K} \right) \right).
\end{eqnarray*}
\end{lem}

 {\bf Proof.} We use a "peeling device". Let $x>0$. For any integer $j$, we consider
$$\cF_j=\left\{f\in\cF:jx\leq A(f)-A^*<(j+1)x \right\}.$$ Define the empirical process
$$Z_x(f)=\frac{A(f)-A_n(f)-(A(f^*)-A_n(f^*))}{A(f)-A^*+x}.$$
Using Bernstein's inequality and margin assumption
MA($\kappa,c,\cF_0$) to upper bound the variance term, we have
\begin{eqnarray*}
\lefteqn{\mathbb{P}\left[\max_{f\in\cF} Z_x(f)>t \right] \leq
\sum_{j=0}^{+\infty}\mathbb{P}\left[\max_{f\in\cF_j}
Z_x(f)>t \right]} \\
& \leq & \sum_{j=0}^{+\infty}
\mathbb{P}\Big[\max_{f\in\cF_j}A(f)-A_n(f)-(A(f^*)-A_n(f^*))>t(j+1)x
\Big]\\ & \leq &
M\sum_{j=0}^{+\infty}\exp\Big(-\frac{n[t(j+1)x]^2}{2c((j+1)x)^{1/\kappa}+(2K/3)t(j+1)x}
\Big)\\ &\leq &
M\Big(\sum_{j=0}^{+\infty}\exp\Big(-\frac{n(tx)^2(j+1)^{2-1/\kappa}}{4cx^{1/\kappa}}
\Big) +\exp\Big(-(j+1)\frac{3ntx}{4K} \Big) \Big)\\
& \leq &
M\Big(\exp\left(-\frac{nt^2x^{2-1/\kappa}}{4c}\right)+\exp\left(-\frac{3ntx}{4K}
\right)\Big)\\ & &
+M\int_{1}^{+\infty}\Big(\exp\left(-\frac{nt^2x^{2-1/\kappa}}{4c}u^{2-1/\kappa}
\right)+\exp\left(-\frac{3ntx}{4K}u \right)\Big)du.\\
\end{eqnarray*}
Lemma \ref{intalpha} completes the proof.

\begin{lem}\label{intalpha}
Let $\alpha\geq1$ and $a,b>0$. An integration by part yields
$$\int_a^{+\infty}\exp\left(-bt^\alpha \right)dt\leq \frac{\exp(-ba^\alpha)}{\alpha b
a^{\alpha-1}}$$
\end{lem}

{\bf{Proof of Corollaries \ref{oraclerégression} and
\ref{oracledensité}.}} In the bounded regression setup, any
probability distribution $\pi$ on $\cX\times[0,1]$ satisfies the
margin assumption MA($1,16,\cF_1$), where $\cF_1$ is the set of
all measurable functions from $\cX$ to $[0,1]$. In density
estimation with the integrated squared risk, any probability
measure $\pi$ on $(\cZ,\cT)$, absolutely continuous w.r.t. the
measure $\mu$ with one version of its density a.s. bounded by a
constant $B\geq1$, satisfies the margin assumption
MA($1,16B^2,\cF_B$) where $\cF_B$ is the set of all non-negative
function $f\in L^2(\cZ,\cT,\mu)$ bounded by $B$. To complete the
proof we use that for any $\epsilon>0$,
$$\Big(\frac{\cB(\cF_0,\pi,Q)\log
M}{\beta_1n}\Big)^{1/2}\leq \epsilon \cB(\cF_0,\pi,Q)+\frac{\log
M}{\beta_2n\epsilon}$$ and in both cases $f\longmapsto Q(z,f)$ is
convex for any $z\in\cZ$.

{\bf{Proof of Theorem \ref{multiperf}.}} We apply Theorem
\ref{oracledensité}, with $\epsilon=1$, to the multi-thresholding
estimator $\hat f_n$ defined in (\ref{multi-thresholding}). Since
the density function $f^*$ to estimate takes its values in
$[0,B]$, $Card(\Lambda_n)= \log n$ and $m\geq n/2$, we have,
conditionally to the first subsample $D_m$,
\begin{eqnarray*}
\lefteqn{\mathbb{E}[\|f^*-\hat{f}_n\|_{L^2([0,1])}^{2} \ | D_m ]} \\
& \leq& 2
\min_{u\in \Lambda_n}(||f^*-
h_{0,B}(\hat f_{v_u}(D_m,.))||_{L^2([0,1])}^{2})+\frac{4(\log n)\log
(\log
n)}{\beta_2n}\\
&\leq & 2 \min_{u\in \Lambda_n}(||f^*-\hat
f_{v_u}(D_m,.)||_{L^2([0,1])}^{2})+\frac{4(\log n)\log (\log
n)}{\beta_2n},\end{eqnarray*} where $h_{0,B}$ is the projection
function introduced in (\ref{EquaProjFunc}) and $\beta_2$ is given
in (\ref{Beta2}). Now, for any $s>0$, let us consider $j_s$ an
integer in $\Lambda_n$ such that $n^{1/(1+2s)}\le 2^{j_s} <
2n^{1/(1+2s)}$. Since the estimators $\hat \alpha_{j,k}$ and $\hat
\beta_{j,k}$ defined by (\ref{est1}) satisfy the inequalities
(\ref{moment}) and (\ref{deviation}), Theorem \ref{guy} implies
that, for any $p\in [1,\infty]$, $s\in (1/p, N]$, $q\in
[1,\infty]$ and $n$ large enough, we have
\begin{eqnarray*}
\lefteqn{
\sup_{f^*\in {{B}}^{s}_{p,q} (L)}  \mathbb{E}[ \Vert
\tilde{f}  -f^*\Vert^{2}_{L^2([0,1])}] =\sup_{f^*\in {{B}}^{s}_{p,q} (L)}  \mathbb{E}[\mathbb{E}[ \Vert
\tilde{f}  -f^*\Vert^{2}_{L^2([0,1])} \ | D_m ] ] } \\ & \leq&   2 \sup_{f^*\in
{{B}}^{s}_{p,q} (L)} \mathbb{E}[\min_{u\in \Lambda_n}(||f^*-\hat
f_{v_u}(D_m,.)||^{2}_{L^2([0,1])}]+\frac{4(\log n)\log
(\log n)}{\beta_2n}\\
& \leq & 2 \sup_{f^*\in {{B}}^{s}_{p,q}(L)}
\mathbb{E}[||f^*-\hat
f_{v_{j_s}}(D_m,.)||^{2}_{L^2([0,1])}]+\frac{4(\log n)\log
(\log n)}{\beta_2n}\\
& \leq & C n^{-2s/(1+2s)}.
\end{eqnarray*}
This completes the proof of Theorem \ref{multiperf}.

{\bf{Proof of Theorem \ref{multiperf2}}.} The proof of Theorem
\ref{multiperf2} is similar to the proof of Theorem
\ref{multiperf}. We only need to prove that, for any
$j\in \{\tau,...,j_1\}$ and $k\in \{0,...,2^j-1\}$, the estimators $\hat
\alpha_{j,k}$ and $\hat \beta_{j,k}$ defined by (\ref{est2})
satisfy the inequalities (\ref{moment}) and (\ref{deviation}).
First of all, let us notice that the random variables
$Y_1\psi_{j,k}(X_1) ,...,Y_n\psi_{j,k}(X_n)$ are i.i.d and that
there $m-$th moment, for $m\geq 2$, satisfies
$$\mathbb{E}(|\psi_{j,k}(X_1)|^m) \le \Vert
\psi\Vert^{m-2}_{\infty}2^{j(m/2-1)}
\mathbb{E}(|\psi_{j,k}(X_1)|^2)=\Vert
\psi\Vert^{m-2}_{\infty}2^{j(m/2-1)}.$$

For the first inequality (cf. inequality (\ref{moment})),
Rosenthal's inequality (see \cite[p.241]{hardle}) yields, for any
$j\in \{\tau,...,j_1\}$,
\begin{eqnarray*}
\mathbb{E}(|\hat \beta_{j,k}-\beta_{j,k}|^4)& \le & C (n^{-3}\mathbb{E}(|Y_1\psi_{j,k}(X_1)|^4) +n^{-2}[\mathbb{E}(|Y_1\psi_{j,k}(X_1)|^2)]^2 )\\
& \le & C \Vert Y\Vert_{\infty}^4\Vert
\psi\Vert_{\infty}^4(n^{-3}2^{j_1}+n^{-2})\le C n^{-2}.
\end{eqnarray*}

For second inequality (cf. inequality (\ref{deviation})),
Bernstein's inequality yields
\begin{equation*}
\mathbb{P}\Big(2\sqrt{n}|\hat
\beta_{j,k}-\beta_{j,k}|\geq\rho\sqrt{a}\Big)\leq
 2\exp\Big(- \frac{\rho^2a}{8\sigma^2+(8/3)M\rho\sqrt{a}/(2\sqrt{n})}\Big),
\end{equation*}
where $a\in \{\tau,...,j_1\}$, $\rho\in
(0,\infty)$,
\begin{eqnarray*}
M & = & \Vert Y\psi_{j,k}(X)-\beta_{j,k} \Vert_{\infty}\le 2^{j/2}
\Vert Y\Vert_{\infty} \Vert\psi\Vert_{\infty}+\Vert f^*\Vert_{L^2([0,1])}^2\\
&\le
&2^{j_1/2}(\Vert \psi\Vert_{\infty}+1) \le  2^{1/2}(n/\log n)^{1/2} (\Vert \psi\Vert_{\infty}+1),
\end{eqnarray*}and
$$\sigma^2=\mathbb{E}(|Y_1\psi_{j,k}(X_1)-\beta_{j,k}|^2) \leq
 \mathbb{E}(|Y_1\psi_{j,k}(X_1)|^2)\le \Vert Y\Vert^2_{\infty}\le 1. $$
Since $a\le \log n$, we complete the proof by seeing that for $\rho$ large enough, we
have
\begin{equation*}
\exp\Big(-\frac{\rho^2a}{8\sigma^2+(8/3)M\rho\sqrt{a}/(2\sqrt{n})}\Big)
\leq 2^{-4a}.
\end{equation*}

\bibliographystyle{plain}
\bibliography{biblioclassification}

\end{document}